\documentclass{birkjour}
\usepackage{amsfonts,amsmath,upref,amssymb,amsthm,amsxtra,amscd,enumerate}

\DeclareMathOperator{\diver}{div}
\DeclareMathOperator{\curl}{curl}

\DeclareMathOperator{\supp}{supp}

\newcommand{\dx}{\,\mathrm{d}x}

\newcommand{\dt}{\,\mathrm{d}t}

\newcommand{\es}{\mathcal{S}}

\newcommand{\ese}{\mathcal{S}^{\varepsilon}}
\newcommand{\ue}{u^{\varepsilon}}
\newcommand{\mue}{\mu^\varepsilon}
\newcommand{\Phie}{\Phi^{\varepsilon}}
\newcommand{\pie}{\pi^{\varepsilon}}

\newcommand{\pt}{\partial_\tau}
\newcommand{\pn}{\partial_\nu}
\newcommand{\pat}{\partial_t}
\newcommand{\esssup}{\operatorname{esssup}}

\newtheorem{theorem}{Theorem}[section]
\newtheorem{lemma}[theorem]{Lemma}
\newtheorem{definition}[theorem]{Definition}

\newtheorem{remark}[theorem]{Remark}

\numberwithin{equation}{section}

\begin{document}
\title[H\"older continuity for shear-thinning fluids]{H\"older continuity of velocity gradients for shear-thinning fluids under perfect slip boundary conditions}
\date{\today}
\author{V\' aclav M\' acha}\address{Industry-University Research Center,\br Yonsei University 50 Yonsei-ro Seodaemun-gu,\br Seoul, 03722, Republic of Korea}\email{macha@math.cas.cz}
\thanks{The research of V\'aclav M\'acha has been supported by the grant NRF-20151009350.}
\author{Jakub Tich\'y}\address{Czech Technical University in Prague,\br Faculty of Information Technology,\br Th\'akurova 9, 160 00, Praha 6, Czech Republic}\email{jakub.tichy@fit.cvut.cz}


\begin{abstract} This paper is concerned with non-stationary flows of shear-thinning fluids in a bounded two-dimensional $\mathcal C^{2,1}$ domain. Assuming perfect slip boundary conditions, we provide a proof of the existence of a solution with the H\"older continuous velocity gradients and pressure under condition that a stress tensor satisfies power-law with growth $p\in[5/3;2]$.
\end{abstract}

\maketitle

\section{Introduction}

The main issue of this paper is to provide a proof of the regularity properties of weak solutions to the system of nonlinear partial differential equations describing evolutionary planar flows of incompressible shear-thinning fluids. The flow of such fluid is governed by the following system:
\begin{equation}\label{ns}
\begin{split}
\partial_t u -\diver \es(Du) + (u\cdot\nabla)u + \nabla \pi &= f,\quad \diver u=0\,\,\textrm{in}\,\,Q,\\
u(0,\cdot)&=u_0\,\, \textrm{in}\,\,\Omega.
\end{split}
\end{equation}
Unknowns $u$ and $\pi$ stand for a velocity and a pressure respectively. Further, $f$ is a density of volume forces and $\es(Du)$ denotes an extra stress tensor. $Du$ is a symmetric part of the velocity gradient, i.e. $Du=\frac{1}{2}[\nabla u + (\nabla u)^\top]$. A time space cylinder is denoted by $Q=I\times\Omega$, where $\Omega \subset \mathbb{R}^2$ is a bounded domain with $\mathcal C^{2,1}$ boundary and $I=(0,T)$ denotes a finite time interval.
We complement \eqref{ns} with the perfect slip boundary conditions
\begin{equation}\label{bcps}
u\cdot\nu=0,\quad [\es(Du)\nu]\cdot\tau=0\quad \textrm{on}\,I\times\partial\Omega,
\end{equation}
which can be considered as a limit case of Navier's slip boundary conditions:
$$
u\cdot\nu=0,\quad \alpha [\es(Du)\nu]\cdot\tau + (1-\alpha)u_\tau = 0\quad \alpha \in[0,1] \textrm{ on }\partial \Omega,
$$
when $\alpha = 1$. Homogeneous Dirichlet boundary conditions, the most often studied case, can be obtained from Navier's slip boundary conditinons when $\alpha = 0$. Here $\tau$ is the tangent vector and $\nu$ is the outward normal to $\partial\Omega$.

The constitutive relation for the extra stress tensor $\es$ is described with the help of generalized viscosity $\mu: [0,\infty) \rightarrow [0,\infty)$ and is of the form 
$$
\es(Du):=\mu(|Du|)Du.
$$ 
Moreover, $\es$ is assumed to possess $p-$potential structure with $p\in(1,2]$, i.e. we assume the existence of the scalar potential $\Phi:[0,\infty) \mapsto [0,\infty)$ to the stress tensor $\mathcal S$ that can be constructed:
\begin{equation*}
\mathcal S(A)=\partial_{A}\Phi(|A|)= \Phi '(|A|)\frac{A}{|A|} \quad \forall A\in\mathbb{R}^{2\times 2}_{sym},
\end{equation*}
such that $\Phi \in \mathcal{C}^{1,1}((0,\infty))\cap\mathcal{C}^{1}([0,\infty))$, $\Phi(0)=0$ and there exist $p\in(1,2]$ and $0<c_1\le c_2$ such that for all $A,B\in\mathbb{R}^{2\times 2}_{sym}$
\begin{equation}\label{as}
c_1(1+|A|^2)^{\frac{p-2}{2}}|B|^2 \le \partial^2_A\Phi(|A|):B\otimes B \le c_2 (1+|A|^2)^{\frac{p-2}{2}}|B|^2.
\end{equation}

H\"older continuity of velocity gradients is an important issue in the mathematical theory of fluid mechanics. Several fields of application of this property of solutions can be considered, for instance exponential attractors. With this kind of regularity it is possible to show the differentiability of the solution operator with respect to the initial condition, which is the key technical step in the method of Lyapunov exponents. Differentiability of the solution is equivalent to the linearisation  of the equation around particular solution which is used to study infinitesimal volume elements and leads to sharp dimension estimates of the global attractor. For more details c.f. \cite{kp}.

Optimal control theory is an another example of possible application. Global regularity results that guarantee boundedness of velocity gradients are needed in order to establish the existence of the weak solution for adjoint equation to the original problem and for linearised models. These results are closely related to the regularity of coefficients in the main part of the associated differential operators and enables to derive corresponding optimality conditions, as is done for example in \cite{wr}. We refer also to \cite{arada1} and \cite{arada2} as an example 
of the work where author is dealing with the lack of the regularity result.

Among many results regarding the regularity theory for non-Newtonian fluids we would like to accentuate those which are dealing with the H\"older continuity of velocity gradients in a bounded domain $\Omega \subset \mathbb R^2$. Starting with the stationary variant of the system \eqref{ns} there are two papers by the authors Kaplick\'y, M\'alek and Star\'a where the H\"older continuity of the velocity gradients is shown for different kinds of boundary conditions and for the case $p<2$. In \cite{kms2} the equation of motion is equipped with the homogeneous Dirichlet conditions and the article \cite{kms1} concerns the non-homogeneous Dirichlet boundary conditions with two types of restriction on boundary data and perfect slip boundary conditions.
     
     In the evolutionary case the first successful attempt to obtain $\mathcal C^{1,\alpha}$ solutions which were global in time was presented in \cite{ser} by Seregin . It was done for $p=2$ only and the boundedness of second and third derivatives of the potential $\Phi$ was required. Under the same assumption the result was extended in \cite{ls}, where the Lipschitz continuity of $\nabla u$ was obtained. Kaplick\'y, M\'alek and Star\'a showed in \cite{kms3} the H\"older continuity of velocity gradients for periodic boundary conditions provided $p\in(4/3,2)$. Later Kaplick\'y obtained in \cite{kap} similar results for the system \eqref{ns} equipped with the homogeneous Dirichlet boundary conditions and $p\in[2,4)$. Using the same structure of the proof the results were extended in \cite{tic} to perfect slip boundary conditions and $p\in[2,\infty)$. By the knowledge of the authors the case $p<2$ is not covered for the perfect slip boundary conditions. Moreover, the shear-thinning fluids, which corresponds to the case when the the generalized viscosity is decreasing function, i. e. $p<2$, are from the view of application more interesting than shear-thickening fluids, i. e. $p>2$. This is the reason we pay attention for this case.

\subsection{Notation}

In this paper we use a standard notation for usual Lebesgue spaces $L^q(\Omega)$, Sobolev-Slobodecki\u\i $\,$ spaces $W^{s,q}(\Omega)$, Besov spaces $B^s_{p,q}(\Omega)$, Bochner spaces $L^q(I,X)$ and $W^{\alpha,q}(I,X(\Omega))$, where $X$ is a Banach space and $\alpha \in (0,1)$, $p,q\in[1,\infty)$, $s \in\mathbb{R}$. Their norms are denoted as $\|\cdot\|_q$, $\|\cdot\|_{s,q}$, $\|\cdot\|_{L^q(X)}$ and $\|\cdot\|_{W^{\alpha,q}(X)}$ respectively. We abbreviate the notation by omitting $\Omega$ and $I$ in case of function spaces, i.e. for instance $W^{s,q}:=W^{s,q}(\Omega)$ and $W^{\alpha,q}(X):= W^{\alpha,q}(I,X(\Omega))$. We do not use different notation for space of scalar, vector-valued or tensor valued functions as far as there will be no misunderstandings. We also write $\int_\Omega f$ instead of $\int_\Omega f(x)\dx$. The symbol $\mathcal C^{0,\alpha}(Q)$, $\alpha\in(0,1)$ stands for a usual space of H\"older continuous functions and $BUC$ means bounded and uniformly continuous functions. 

Let $\mathcal{C}^{\infty}_{\sigma}=\{\varphi \in\mathcal{C}^{\infty},\,\diver \varphi=0\,\textrm{in}\,\Omega,\,\varphi\cdot\nu=0\,\textrm{on}\,\partial\Omega \}$. 
Since the domain $\Omega$ is in our case at least $\mathcal{C}^{2,1}$, we can define $L^q_\sigma$ and $W^{1,q}_{\sigma}$ as follows:

\begin{displaymath}
\begin{split}
L^q_\sigma&=\{\varphi\in L^q,\,\diver \varphi =0\,\textrm{in}\,\Omega,\,\varphi\cdot\nu=0\,\textrm{on}\,\partial\Omega\},\\
W^{1,q}_{\sigma}&=\{\varphi \in W^{1,q},\,\diver\varphi=0,\,\textrm{in}\,\Omega,\,\,\varphi\cdot \nu=0\, \textrm{on}\,\partial\Omega\}.
\end{split}
\end{displaymath}
By $\langle\cdot,\cdot\rangle$ we denote the duality pairing between Banach space $X$ and its dual $X'$.
Set $W^{-1,p'}_\sigma:=(W^{1,p}_\sigma)'$.

\subsection{Main result}

We begin with the definition of the weak solution to the problem (\ref{ns}) with (\ref{bcps}).

\begin{definition}
Let $f \in L^{p'}(I,W^{-1,p'}_\sigma)$, $p\in(1,2]$ and $u_0 \in L^2$. We say that a function $u: Q \mapsto \mathbb{R}^2$ is a weak solution to the problem \eqref{ns} with \eqref{bcps}, if $u \in L^\infty(L^2)\cap L^p(W^{1,p}_\sigma)$, $\partial_t u \in L^{p'}(W^{-1,p'}_\sigma)$, $u(0,\cdot)=u_0$ in $L^2$ and weak formulation
\begin{equation*}
\int_I\langle\partial_t u,\varphi\rangle + \int_Q\es(Du)\!:\!D\varphi + \int_Q(u \otimes u)\nabla\varphi - \int_Q \pi \diver \varphi = \int_I \langle f,\varphi\rangle
\end{equation*}
holds for all $\varphi \in \mathcal C^\infty_0(-\infty,T,\mathcal C^\infty)$.
\end{definition}

It is well known that $u\in \mathcal C(\bar{I},L^2)$ and the initial value problem is well posed. Since $p<2$, we can't take test functions from the space from $L^p(I,W^{1,p}_\sigma)$ and uniqueness of the weak solution is not as obvious as for the case $p>2$ and needs to be shown. The existence theory for periodic boundary conditions can be found for instance in \cite[Chapter 5]{mnrr}.

The main result of this paper is formulated as follows.

\begin{theorem}\label{thm1}
Let $\Omega \subset \mathbb{R}^2$ be a bounded non-circular $\mathcal{C}^{2,1}$ domain, $p\in[5/3,2]$ and \eqref{as} holds. Let $u_0 \in W^{2+\beta,2}$ for $\beta\in(0,1/4)$, $\diver u_0=0$, $f\in L^{\infty}(L^{q_0})$ and $\partial_t f \in L^2(L^{p'})\cap L^{q_0}(W^{-1,{q_0'}}_{\sigma})$ for some $q_0>2$. 
Then there exists a solution $(u,\pi)$ of \eqref{ns} with \eqref{bcps}, such that for some $\alpha>0$
$$
\nabla u,\pi \in \mathcal C^{0,\alpha}(Q).
$$
\end{theorem}

Theorem \ref{thm1} is proven in \cite[Section 4]{tic} provided $p=2$. A strategy of the proof of Theorem \ref{thm1} for $p\in[5/3,2)$ is following. We approximate the stress tensor $\es$ by truncation of the generalized viscosity $\mu$ and consider \eqref{ns} with this approximated tensor $\ese$ instead of the original one $\es$. This is done in Section \ref{app.sys}. In Section~\ref{unif.est} we obtain estimates independent of the approximation parameter $\varepsilon$. Further, we improve time regularity in Section~\ref{imp.tim.reg} by application of $L^p$ theory result for linear Stokes system and after subsuming time derivative of solution into the right hand side we enhance space regularity in Section~\ref{imp.spa.reg} using $L^p$ theory for stationary Stokes system. In the last section we explain how to pass from approximate system to the original problem.

\subsection{Used inequalities}

We recall several inequalities which are commonly used throughout this paper, often without specific remark.
\begin{itemize}
\item Korn's inequality (c. f. \cite[Lemma 2.7]{KaTi}): For a non-circular $\mathcal C^{0,1}$ domain $\Omega\subset \mathbb R^2$ and arbitrary $q\in (1,\infty)$ there exists $c>0$ such that for all $v\in W^{1,q}$, $v\cdot \nu = 0$
$$
\|\nabla v\|_p\leq c\|D v\|_p.
$$
\item Replacement of the full second gradient by a sum of three gradients of the symmetric gradient:
It can be verified by a simple calculation that for arbitrary function $v$  it holds
$$
|\nabla^2 v| \leq 3|\nabla D v|.
$$

\item Poincar\'e's inequality: For a bounded $\mathcal C^{0,1}$ domain $\Omega$ and arbitrary $q\in (1,\infty)$ there exists $c>0$ such that
\begin{equation*}
\begin{split}
\|v\|_{1,q}&\leq c\|\nabla v \|_q \qquad \mbox{for}\ v\in W^{1,q},\ v\cdot \nu = 0 \textrm{ on }\partial\Omega,\\
\|v\|_{2,q}&\leq c\|\nabla^2 v\|_q \qquad \mbox{for}\ v\in W^{2,q},\ v\cdot \nu = 0 \textrm{ on }\partial\Omega.\\
\end{split}
\end{equation*}
These inequalities can be proven by a simple contradiction argument.
\item Gagliardo-Nirenberg's inequality (c. f. \cite[Lecture II]{niren}): Let $\Omega \subset \mathbb R^2$, $k\in \mathbb N$, $\beta\in \mathbb N\cup \{0\}$ and let $r,\ s,\ q \in (1,\infty)$ satisfy
$$
\frac 1s = \frac \beta 2 + \lambda\left(\frac 1q - \frac k2\right) + (1-\lambda)\frac 1r,\quad \frac \beta k \leq \lambda \leq 1,\quad \beta\leq k-1.
$$
Then there exists $c$ such that
$$
\|\nabla^\beta v\|_{s}\leq c \|v\|_{k,q}^\lambda \|v\|_{r}^{1-\lambda},\quad  \mbox{for all}\ v\in W^{k,q},
$$
provided $k-\beta-\frac 2q$ is not a negative integer.
\end{itemize}

\section{Approximate system}\label{app.sys}
Instead of \eqref{ns} we consider the following approximate system

\begin{equation}\label{ap.sys}
\begin{array}{rcl}
\pat \ue - \diver \ese(D\ue) + (\ue\cdot\nabla)\ue + \nabla \pie &=& f\mbox{ on }Q,\\
\diver \ue &=& 0\mbox{ on }Q,\\
\ue(0,\cdot) &=& u_0 \mbox{ on }\Omega,
\end{array}
\end{equation}
together with perfect slip boundary conditions \eqref{bcps}. The stress tensor 
$$\ese(D\ue) = \mue(|D\ue|)D\ue$$ 
is defined via
$$
\mue(s) = \max \{\varepsilon, \mu(s)\},\quad s,\varepsilon>0
$$
and the scalar potential $\Phie$ to $\ese$ is defined by $\Phie(s):=\int_0^s\mue(t)t\dt$. Due to the truncation of the generalized viscosity the growth conditions \eqref{as} can be rewritten as
\begin{equation}\label{trunc.growth.con}
\gamma_1(\varepsilon) |B|^2\leq \partial_A^2 \Phie(|A|) : B\otimes B \leq \gamma_2|B|^2.
\end{equation}
Here $\gamma_2$ is independent of $\varepsilon$ and $\gamma_1$ can be determined as follows
$$
\gamma_1(\varepsilon) =  \max \left\{\varepsilon, (1+|A|^2)^{\frac{p-2}{2}}\right\}.
$$

We recall an estimate obtained in \cite[Section 4]{tic} in proof of the regularity result to equations with quadratic growth:
\begin{equation}\label{first.apriori}
\sup_{t\in I}\|\ue(t,\cdot)\|_2 \leq c \|f\|_{-1,p},
\end{equation}
where $c$ is $\varepsilon$-independent.

Note that for every $\varepsilon \in (0,1)$ gradient of the solution $\ue$ to \eqref{ap.sys} is H\"older continuous. Thus, all the calculations in the following sections make sense. Our goal is to derive $\varepsilon$-independent estimates.

\section{Boundary description}\label{boundary.desc}
In order to discuss boundary regularity in the following steps, we need a suitable description of the boundary $\partial\Omega$. Let us denote $x=(x_1,x_2)$. We suppose that $\Omega \in \mathcal{C}^{2,1}$, therefore there exists $c_0>0$ such that for all $a_0>0$ there exists $n_0$ points $P\in \partial\Omega$, $r>0$ and open smooth set $\Omega_0 \subset \subset \Omega$ that we have
\begin{equation}\label{bodP}
\Omega \subset \Omega_0 \cup \bigcup_P B_{r}(P)
\end{equation}
and for each point $P\in\partial \Omega$ there exists local system of coordinates for which $P=0$ and the boundary $\partial \Omega$ is locally described by $\mathcal{C}^{2,1}$ mapping $a_P$ that for $x_1 \in(-3r,3r)$ fulfills
\begin{displaymath}
\begin{split}
&x \in \partial \Omega \Leftrightarrow x_2=a_P(x_1), \\ 
&B_{3r}(P)\cap\Omega = \{x\in B_{3r}(P)\,\textrm{and}\,x_2>a_P(x_1)\}=:\Omega_{3r}^P,
\\
&\partial_1 a_P(0)=0, \quad |\partial_1 a_P(x_1)|\le a_0,\quad |\partial_1^2 a_P(x_1)|+|\partial_1^3 a_P(x_1)|\le c_0.
\end{split}
\end{displaymath}
Points $P$ can be divided into $k$ groups such that in each group $\Omega_{3r}^P$ are disjoint and $k$ depends only on dimension $n$. Let the cut-off function $\xi_P(x) \in \mathcal{C}^{\infty}(B_{3r}(P))$ and reaches values
\begin{displaymath}
\xi_P(x) \left\{ \begin{array}{ll}
=1  &x \in B_r(P), \\
\in (0,1) &x \in B_{2r}(P)\setminus B_r(P),\quad  \\
=0 &x\in \mathbb{R}^2\setminus B_{2r}(P).
\end{array} \right.
\end{displaymath}
Next, we assume that we work in the coordinate system corresponding to $P$. Particularly, $P=0$. Let us fix $P$ and drop for simplicity the index $P$. The tangent vector and the outer normal vector to $\partial\Omega$ are defined as
$$
\tau=\big(1,\partial_1 a(x_1)\big), \quad \nu=\big(\partial_1a(x_1),-1\big),
$$
tangent and normal derivatives as
$$
\pt = \partial_1 + \partial_1 a(x_1)\partial_2, \quad \pn = -\partial_2 + \partial_1 a(x_1)\partial_1.
$$

Because the function $a(x_1)$ doesn't depend on $x_2$, we will abbreviate the notation as follows $a':=\partial_1a(x_1)$ and $a'':=\partial_1^2a(x_1)$.

\section{Uniform estimates}
\label{unif.est}
The goal of this section is to derive space and time regularity results by estimates uniform with respect to $\varepsilon$. These kind of estimates are commonly proven only locally. We omit interior estimates since that can be done similarly as in \cite[Section~4]{kms3}, there we would test the equation of motion by $\Delta\ue$, more precisely by $\nabla^\bot(\xi^2 \curl \ue)$. We focus only on the estimates near boundary, therefore we need a test function appropriate for perfect slip boundary conditions. To estimate all lower order terms we collect all local estimates and continue working on the whole domain $\Omega$.

For a scalar function $g$ and a vector function $v$ we denote 
\begin{eqnarray*}
\nabla^\bot g=(\partial_2 g, -\partial_1 g),\qquad \curl v = \partial_1 v_2 - \partial_2 v_1.
\end{eqnarray*}

As a test function in the weak formulation of \eqref{ap.sys} we take
\begin{equation}\label{test.dd}
\varphi= \nabla^\bot(\xi^2\Theta), 
\end{equation}
where 
\begin{equation}\label{def.theta}
\begin{split}
\Theta:&= (\pn\ue\cdot\tau -\ue\cdot\pt\nu) -\pt(\ue\cdot\nu)\\
&= (1+a'^2)\curl \ue - 2a'' \ue_1.
\end{split}
\end{equation}

Note that test function $\varphi$ in \eqref{test.dd} fulfills $\diver \varphi = 0$, because for every scalar $g$ smooth enough it holds $\diver \nabla^\bot g = 0$. Second inportant property of $\varphi$ enables to integrate by parts without any boundary integrals. Indeed, as one can easily see, it holds $\Theta = 0$ on $\partial\Omega$ (c.f. \cite[Proposition~4.4]{KaTi}).

\subsection{Local estimates}

We multiply \eqref{ap.sys} by $\varphi$ defined via \eqref{test.dd} and we integrate over $(0,\tau)\subset I$ in order to get
\begin{equation}\label{flat.test}
\begin{split}
\int_{(0,\tau)\times\Omega_{3r}}  \partial_t \ue \nabla^\bot(\xi^2 \Theta) - \int_{(0,\tau)\times\Omega_{3r}} \diver \ese(D\ue)\nabla^\bot(\xi^2 \Theta) \\+ \int_{(0,\tau)\times\Omega_{3r}} (\ue\cdot\nabla)\ue \nabla^\bot(\xi^2 \Theta) = \int_{(0,\tau)\times\Omega_{3r}} f \cdot \nabla^T(\xi^2 \Theta).
\end{split}
\end{equation}

The term containing the pressure vanished, because the test function was constructed in order to fulfil $\diver \varphi = 0$. We estimate all terms of \eqref{flat.test} in the similar way like in \cite[Proof of Theorem 1.8]{KaTi}, where authors studied stationary Stokes problem. Therefore we should focus especially on the term containing time derivative and convective term. Estimates of the elliptic term, i.e.  the second term of \eqref{flat.test}, demands more difficult calculations, although the method is the same as method presented here. Since it is done in \cite[Section~4]{KaTi}, we omit it.  The desired term obtained from the elliptic term is $\int_{(0,\tau)\times\Omega_{3r}}\mue(|D\ue|)|\nabla D\ue|^2 \xi^2$. Except terms controlled by apriori estimates there appear lower order term which could be estimated by $c\int_{(0,\tau)\times\Omega_{3r}}\mue(|D\ue|)|\nabla D\ue|(|\nabla \ue| + |\ue|)\xi$ and small same order terms that can be subsumed into $\int_{(0,\tau)\times\Omega_{3r}}\mue(|D\ue|)|\nabla D\ue|^2 \xi^2$ due to the smallness provided by the suitable choice of the function $a$ locally describing the boundary $\partial\Omega$.

We use integration by parts in the first term of \eqref{flat.test} in order to get
\begin{multline}\label{ieoasg1}
\int_{(0,\tau)\times \Omega_{3r}}\pat \ue \nabla^\bot (\xi^2 \Theta) = \int_{(0,\tau)\times \Omega_{3r}} \pat \curl \ue [(1+a'^2)\curl \ue - 2a'' \ue_1]\xi^2\\
 = \int_0^\tau \pat \int_{\Omega_{3r}} \frac 12 (1+a'^2) |\curl \ue|^2\xi^2 - \int_{(0,\tau)\times \Omega_{3r}} 2a'' \pat(\curl \ue)\ue_1\xi^2\\
= \int_{\Omega_{3r}} \frac 12 (1+a'^2) |\curl \ue(t)|^2\xi^2 - \int_{\Omega_{3r}}\frac 12 (1+a'^2) |\curl \ue(0)|^2\xi^2\\
 -\int_{\Omega_{3r}} 2a'' [\curl \ue(t)] \ue_1(t) \xi^2 + \int_{\Omega_{3r}} 2a'' [\curl \ue(0)]\ue_1(0) \xi^2\\ +\int_{(0,\tau)\times \Omega_{3r}} 2a'' \curl \ue \pat \ue_1 \xi^2 = \sum_{k=1}^5 \mathcal I_k.
\end{multline}
The term $\mathcal I_1$ has a suitable sign, $\mathcal I_2$ can be handled by initial conditions as well as the term $\mathcal I_4$. Further, we apply Young's inequality to $\mathcal I_3$ and we arrive at
$$
|\mathcal I_3| =\left|\int_{\Omega_{3r}} 2a'' [\curl \ue(t)]\ue_1(t)\xi^2\right|\leq \delta \int_{\Omega_{3r}}|\curl \ue(t)|^2\xi^2 + c_\delta \int_{\Omega_{3r}}|\ue_1(t)|^2\xi^2,
$$
where the first term on the right hand side can be absorbed by $\mathcal I_1$ and the second term can be estimated by \eqref{first.apriori}.

It remains to estimate $\mathcal I_5$. Recalling that $a\in \mathcal C^{2,1}$ we use in the first step the definition of $\Theta$ \eqref{def.theta} to replace $\curl \ue$.

\begin{multline*}
\mathcal I_5 = \int_{(0,\tau)\times \Omega_{3r}} 2a'' \curl \ue\pat\ue_1\xi^2 = \int_{(0,\tau)\times \Omega_{3r}} \frac{2a''}{1+a'^2} \Theta \pat \ue_1\xi^2 \\ + \int_{(0,\tau)\times \Omega_{3r}} \frac{2a''^2}{1+a'^2}\ue_1\pat\ue_1\xi^2 = \mathcal I_6 + \mathcal I_7.
\end{multline*}
The term $\mathcal{I}_7$ is bounded due to appriori estimates \eqref{first.apriori}:
$$
|\mathcal I_7| \leq c (\|\ue(t)\|_{2} + \|\ue(0)\|_2)\leq c.
$$

We denote by $\mathcal P$ the Helmholtz projection, which is a bounded operator on space $L^q$ and $W^{1,q}$ for every $q\in (1,\infty)$ (c.f. \cite{cantor}). We use an abbreviation $H(\Theta) = \mathcal P\left(\frac{2a''}{1+a'^2}\xi^2\Theta,0\right)$. As $\diver \pat \ue=0$, the term $\mathcal I_6$ can be treated as follows:
\begin{multline*}
\mathcal I_6 = \int_{(0,\tau)\times \Omega_{3r}}\pat \ue \mathcal P\left(\frac{2a''}{1+a'^2}\xi^2\Theta,0\right) \\ = \int_{(0,\tau)\times \Omega_{3r}} -\ese(D\ue) \nabla H(\Theta) - (\ue\cdot \nabla)\ue H(\Theta) + f H(\Theta)\\
\leq c\Big(1 + \|\ese(D\ue)\|_{L^{p'}(L^{p'})}^{p'} +  \|\Theta\|_{L^p(0,\tau,W^{1,p}(\Omega_{3r}))}^p\\ + \int_{(0,\tau)\times \Omega_{3r}} |\ue||\nabla \ue||H(\Theta)|\Big).
\end{multline*}

Reminding that $\int_{(0,\tau)\times\Omega_{3r}}|\ese(D\ue)|^{p'}\leq c \|\nabla^2 \ue\|^p_{L^p(L^p)},$ together with $\|\Theta\|_{L^p(W^{1,p})}\leq c\|\nabla Du\|_{L^p(L^p)}$ and $\|H(\Theta)\|_{L^4(\Omega_{3r})}\leq c \|\nabla \ue\|_{L^4(\Omega_{3r})}$ we infer
$$
|\mathcal I_6| \le c\left(1+ \|\nabla^2 \ue\|^p_{L^p(L^p)} + \int_0^\tau \|u\|_2 \|\nabla u\|_4^2\right).
$$

Now we can focus on the convective term, i.e. the third term of \eqref{flat.test}. The aim of this detailed procedure is to show that terms containing second derivatives of $\ue$ did not appear due to the divergence-free condition.
\begin{multline*}
\int_{(0,\tau)\times\Omega_{3r}} (\ue\cdot\nabla)\ue\cdot\varphi\\ = \int_{(0,\tau)\times\Omega_{3r}} \partial_2(\xi^2\Theta)[\ue_1\partial_1\ue_1 + \ue_2\partial_2\ue_1]  - \partial_1(\xi^2\Theta)[\ue_1\partial_1\ue_2 + \ue_2\partial_2\ue_2]  \\
= \int_{(0,\tau)\times\Omega_{3r}} \xi^2\Theta(-\ue_1\partial_2\partial_1\ue_1 - \ue_2\partial_2^2\ue_1 + \ue_1\partial_1^2\ue_2 + \ue_2\partial_1\partial_2\ue_2) =\mathcal J_1,
\end{multline*}
where we used the fact that there arise no boundary terms while integrating by parts since $\Theta = 0$ at $\partial\Omega$. Four terms in $\mathcal J_1$ vanished because $\partial_1\ue_1= -\partial_2\ue_2$. Now we put together terms containing $\ue_1$ and integrate by parts in the direction $x_1$, in other terms we integrate by parts in direction $x_2$. We get
\begin{displaymath}
\begin{split}
\mathcal J_1 = \int_{(0,\tau)\times\Omega_{3r}} \xi^2\Theta (\partial_1\ue_1 \partial_2\ue_1 - \partial_1\ue_1\partial_1\ue_2 + \partial_2\ue_2\partial_2\ue_1 - \partial_2\ue_2\partial_1\ue_2) \\+
\int_{(0,\tau)\times\Omega_{3r}}\partial_1(\xi^2\Theta)[\ue_1\partial_2\ue_1 - \ue_1\partial_1\ue_2] + \partial_2(\xi^2\Theta) [\ue_2\partial_2\ue_1 - \ue_2\partial_1\ue_2]  \\
= \mathcal J_2 - \int_{(0,\tau)\times\Omega_{3r}} [\partial_1(\xi^2\Theta)\ue_1 + \partial_2(\xi^2\Theta)\ue_2]\curl\ue
= \mathcal J_2 + \mathcal J_3.
\end{split}
\end{displaymath}
One can easily see that $\mathcal J_2=0$ since $\diver \ue =0$. Using \eqref{def.theta} we write out $\mathcal J_3$ in the following way
\begin{multline*}
\mathcal J_3 = -\int_{(0,\tau)\times\Omega_{3r}}(1+a'^2)\xi^2[(\partial_1^2\ue_2-\partial_1\partial_2\ue_1)\ue_1 
+(\partial_2\partial_1\ue_2-\partial_2^2\ue_1)\ue_2]\curl\ue\\  
+ \int_{(0,\tau)\times\Omega_{3r}} \big([2a'a''\curl\ue +2a'''\ue_1+2a''\partial_1\ue_1]\ue_1+2a''\partial_2\ue_1\ue_2\big)\xi^2\curl\ue \\
-\int_{(0,\tau)\times\Omega_{3r}} 2\xi\Theta(\partial_1\xi\ue_1 + \partial_2\xi\ue_2)\curl\ue=\mathcal J_4 + \mathcal J_5 + \mathcal J_6.
\end{multline*}

The term $\mathcal J_4$ can be rewritten after integration by parts in $x_1$ direction and $x_2$ direction:
\begin{displaymath}
\begin{split}
\mathcal J_4 = -\frac{1}{2} \int_{(0,\tau)\times\Omega_{3r}}(1+a'^2)\xi^2[u_1\partial_1|\curl\ue|^2 + \ue_2\partial_2|\curl\ue|^2]  \\
=\frac{1}{2}\int_{(0,\tau)\times\Omega_{3r}}(1+a'^2)\xi^2(\partial_1\ue_1 + \partial_2\ue_2)|\curl\ue|^2
\\ +\int_{(0,\tau)\times\Omega_{3r}} [a'a''\xi^2\ue_1 + (1+a'^2)(\ue_1\xi\partial_1\xi + \ue_2\xi\partial_2\xi)]|\curl\ue|^2 \\
-\frac 12 \int_{(0,\tau)}\int_{\partial\Omega \cap \partial \Omega_{3r}} (1+a'^2)\xi^2(\ue_1\nu_1 + \ue_2\nu_2)|\curl\ue|^2= \mathcal J_7 +\mathcal J_{8}+\mathcal J_{9}.
\end{split}
\end{displaymath} 
One can see that $\mathcal J_7=0$, since $\diver \ue=0$ and $\mathcal J_{9}=0$, because $\ue\cdot\nu=0$ at $\partial\Omega$. Thus, collecting previous estimates we arrive at

\begin{multline}\label{con.te}
\left|\int_{(0,\tau)\times\Omega_{3r}} (\ue\cdot\nabla)\ue\cdot\varphi\right| \le |\mathcal J_5| + |\mathcal J_6| + |\mathcal J_{8}| \\\le c\int_{(0,\tau)\times\Omega_{3r}} \big(|\ue||\nabla\ue|^2 + |\ue|^2|\nabla\ue|\big)\xi.
\end{multline}

Summarising calculation of this subsection, i.e. estimates obtained from the elliptic terms, term containing time derivative and convective term, we arrive at
\begin{multline*}
\frac 12\int_{\Omega_{3r}} (1+a'^2)|\curl \ue(\tau)|^2\xi^2 + \int_{(0,\tau)\times\Omega_{3r}} \mue(|D\ue|)|\nabla D\ue|^2\xi^2 \\ \leq c\int_{(0,\tau)\times\Omega_{3r}} |f|(|\nabla^2\ue| + |\nabla\ue|)\xi^2 
+ c\int_{(0,\tau)\times\Omega_{3r}}\mue(|D\ue|)|\nabla D\ue||\nabla \ue|\xi\\ + c\int_{(0,\tau)\times\Omega_{3r}}\big(|\ue||\nabla\ue|^2 +|\ue|^2|\nabla\ue|\big)\xi \\
+ c\left(1 + \|\nabla^2 \ue \|^p_{L^p(L^p)}  + \int_0^\tau \|u\|_2 \|\nabla u\|_4^2 \right).
\end{multline*}

\subsection{Global estimates}

According to Section \ref{boundary.desc} we obtain estimates on the whole domain $\Omega$ by combination of interior regularity result together with estimates on $\Omega_{3r}^P$. In the previous subsection we worked on $\Omega_{3r}^P$ for fixed $P\in\partial\Omega$. We recall that points $P$ are divided into $k$ groups and in each group the sets $\Omega_{3r}^P$ are mutually disjoint. Summing over all points $P$ and setting $\tau = T$ we arrive at

\begin{equation}\label{ieoasg3}
\int_Q \mue (|D\ue|) |\nabla D\ue|^2 \leq c\sum_{i\in \{1,2,3\}} \mathcal M_i +c \|\nabla^2\ue\|_{L^p(Q)}+ c
\end{equation}
where 
\begin{equation*}
\begin{split}
\mathcal M_1 &= \int_Q |f|(|\nabla^2 \ue| + |\ue|),\\
\mathcal M_2 &= \int_Q \mue(|D\ue|) |\nabla D\ue|(|\nabla \ue| + |\ue|),\\
\mathcal M_3 &= \int_Q (|\ue||\nabla \ue|^2 + |\ue|^2|\nabla\ue|) + c\int_0^T\|\ue\|_2\|\nabla \ue\|_4^2.
\end{split}
\end{equation*}

Young's inequality and a bound $\mue(|D\ue|)\leq c$ that holds for every $D\ue$ enable to estimate first two terms as

\begin{equation*}
\begin{split}
\mathcal M_1 &\leq c \left(1 + \|\nabla^2\ue\|^p_{L^p(L^p)}\right),\\
\mathcal M_2 
&\leq \delta \int_Q \mue (|D\ue|)|\nabla D\ue|^2 + c_\delta\int_0^\tau \|\nabla \ue\|_2^2.
\end{split}
\end{equation*}

Further, Gagliardo-Nirenberg's and Poincar\'e's inequlities provides 
\begin{equation}\label{GN}
\|\nabla \ue\|_4 \leq  c\|\nabla^2 \ue\|_p^{\frac{3p}{6p-4}}\|\ue\|_2^{\frac{3p-4}{6p-4}}.
\end{equation}
This can be used in estimate of $\mathcal M_2$ as follows
$$
\int_0^T \|\nabla \ue\|_2^2\leq c\int_0^T \|\nabla \ue\|_4^2\leq c \int_0^T \|\nabla^2\ue\|_{p}^{\frac{3p}{3p-2}}.
$$
But the main task of \eqref{GN} is to help with term $\mathcal M_3$ in the following way 
\begin{multline*}
\mathcal M_3
\leq c\int_0^\tau \|\ue\|_2 \|\nabla \ue\|_4^2 \leq c \int_0^\tau \|\ue\|_2^{\frac{3p-4}{3p-2}+1} \|\nabla^2\ue\|_p^{\frac{3p}{3p-2}} \\ \leq c \int_0^\tau \|\nabla^2\ue\|_p^{\frac{3p}{3p-2}}.
\end{multline*}
The assumption $p\geq \frac 53$ implies $\frac{3p}{3p-2}\leq p$ and thus \eqref{ieoasg3} yields
\begin{equation}\label{ieoasg4}
\int_0^T \mue(|D\ue|)|\nabla D \ue|^2 \leq c + c\|\nabla^2 \ue\|_{L^p(L^p)}^p.
\end{equation}
Since 
\begin{multline}\label{zmnrr}
\int_0^\tau \int_\Omega |\nabla^2 \ue|^p  = \int_0^\tau \int_\Omega (1+|D\ue|^2)^{\frac{(p-2)p}{4}} |\nabla^2 \ue|^p (1+|D\ue|^2)^{\frac{(2-p)p}{4}}\\ \leq \delta \int_0^\tau \int_\Omega \left(1+|D\ue|^2\right)^{\frac{p-2}{2}}|\nabla^2 \ue|^2 + c_\delta \int_0^\tau \int_\Omega \left(1+|D\ue|^2\right)^{\frac p2}\\ \leq \delta\int_0^\tau \int_\Omega \mue (|D\ue|)|\nabla D\ue|^2 + c_\delta,
\end{multline}
we derive from \eqref{ieoasg4} the estimate $\int_0^T \mue(|D\ue|)|\nabla D \ue|^2 \leq c$ and consequently again by \eqref{zmnrr} 
\begin{equation}\label{ieoasgfinal}
\|\nabla^2\ue\|_{L^p(L^p)}\leq c.
\end{equation}


\subsection{Iteration}

We are going to introduce iterative method which improves integrability of $\ue$. Assume that $\|\ue\|_{L^\infty(L^q)}\leq c$ for some $q\in [2,5]$. As \eqref{ieoasgfinal} implies that $\|\ue\|_{L^{\frac 53}(W^{2,\frac 53})}\leq c$, we use Gagliardo-Nirenberg's inequatlity in order to derive
\begin{equation}\label{GN2}
\|\ue\|_{L^{\frac 53(q+1)}(Q)}\leq c.
\end{equation}
Gagliardo-Nirenberg's inequality also gives
\begin{equation}\label{iter1}
\|\nabla \ue\|_{L^{\frac 52}(Q)} \leq c,
\end{equation}
which implies $\ue \cdot \nabla \ue\in L^{\frac{5q+5}{2q+5}}$ uniformly in $\varepsilon$. Further, \eqref{ieoasg4} yields that $\|\diver \ese(D\ue)\|_{L^2(Q)}\leq c$. Since $q\leq 5$, it follows that $\frac{5q+5}{2q+5}\leq 2$ and we derive with help of the Helmholtz decomposition and \eqref{ap.sys} that
$$
\|\pat \ue\|_{L^{\frac{5q+5}{2q+5}}(Q)}\leq c.
$$
For every $t\in (0,T)$  
\begin{multline}\label{methoda}
\|\ue(t)\|_{q+1}^{q+1} = \pat\int_0^t\int_\Omega |\ue|^{q+1} + \|\ue(0)\|_{q + 1}^{q + 1}\\
 \leq c\int_0^t\int_\Omega|\ue|^{q} |\pat \ue|+ \|\ue(0)\|_{q + 1}^{q + 1} \\
\leq \int_0^t\int_\Omega |\ue|^{\frac{5}{3}(q+1)} +\int_0^t\int_\Omega |\pat \ue|^{\frac{5q+5}{2q+5}} + \|\ue(0)\|_{q + 1}^{q + 1} \leq c.
\end{multline}
We have just proved that $\|\ue\|_{L^\infty(L^q)}\leq c \Rightarrow \|\ue\|_{L^\infty(L^{q + 1})}\leq c$ provided $q\in [2,5]$.

Starting with $q=2$ we get $\|\ue\|_{L^\infty(L^{5})}\leq c$ after several iterations. This yields by \eqref{GN2} 
\begin{equation}\label{iter2}
\|\ue\|_{L^{10}(Q)}\leq c.
\end{equation}

\begin{remark}
Instead of the iteration resulting to \eqref{iter2} another approach can be used in order to obtain $\ue \in L^\infty(W^{1,2})$ that would allow us to estimate the convective term in the next section. Note that in \eqref{ieoasg3} there could appear $\|\curl \ue(t,\cdot)\|_2$ on the left hand side yielding $\curl \ue \in L^\infty(L^2)$. Boundedness of $\|\curl \ue\|_2$ together with $\diver \ue=0$ suffices to control $\|\nabla \ue\|_2$. Namely, we could use that $\Delta \ue = \nabla^\bot \curl \ue$ which in the weak formulation looks like
\begin{equation*}\label{remarkweak}
\int_\Omega \nabla \ue \nabla\varphi = \int_\Omega \curl \ue \curl \varphi -\int_{\partial \Omega} \curl \ue \varphi \tau,\quad \forall \varphi \in W^{1,2},\, \varphi\cdot \nu = 0.
\end{equation*}
The right hand side contains boundary integral that did not disappear in case of perfect slip boundary conditions. Since we should have dealt with technical difficulties while estimating the boundary integral or constructing a proper test function we preferred the iteration method. 
\end{remark}

\subsection{Estimates of time derivative}
Having at our disposal \eqref{iter1} and \eqref{iter2} we know that $(\ue \cdot \nabla) \ue\in L^2(Q)$ uniformly in $\varepsilon$, therefore we may multiply \eqref{ap.sys} by $\pat \ue$ and derive 
\begin{equation}\label{casoveodhady}
\begin{split}
\|\pat \ue\|_{L^2(Q)} &\leq c,\\
\|\nabla \ue\|_{L^\infty(L^p)} &\leq c.
\end{split}
\end{equation}

We can also differentiate \eqref{ap.sys} with respect to time, multiply it by $\pat \ue$ and integrate it over $(0,\tau)\times \Omega$ for some $\tau \in [0,T]$. We obtain
\begin{multline}\label{skorokonec}
\frac{1}{2}\int_0^\tau \pat \|\pat \ue\|_2^2 + \int_{(0,\tau)\times \Omega} \mu(|D\ue|) |D\pat \ue|^2 \\ \leq \int_{(0,\tau)\times \Omega} |\pat f\pat\ue| + \int_{(0,\tau)\times \Omega} |\pat \diver(u\otimes u) \pat u|.
\end{multline}
Before estimating right hand side of \eqref{skorokonec}, we modify the convective term
\begin{multline*}
\int_{(0,\tau)\times \Omega} \pat \diver(u\otimes u) \pat u = \int_{(0,\tau)\times \Omega} \pat \ue \nabla \ue \pat \ue + \int_{(0,\tau)\times \Omega} \ue \nabla \frac 12|\pat \ue|^2,
\end{multline*}
where the last term is equal to zero, since after integration by parts we apply $\diver \ue =0$. Thus,
$$
\int_{(0,\tau)\times \Omega} |\pat \diver(u\otimes u) \pat u|= \int_{(0,\tau)\times \Omega} |\pat \ue|^2 |\nabla \ue|\leq \int_0^\tau \|\pat \ue\|_{\frac{20}9}^2 \|\nabla \ue\|_{10}.
$$
Due to the Gagliardo-Nirenberg's inequality $\|\pat \ue\|_{\frac {20}9}\leq c\|\pat \nabla \ue\|_{\frac 53}^{\frac 1{8}}\|\pat \ue\|_2^{\frac {7}{8}}$ and we can use Young's inequality in order to get
\begin{multline*}
\int_0^\tau \|\pat \ue\|_{\frac{20}9}^2 \|\nabla \ue\|_{10} \leq \int_0^\tau \|\pat \nabla \ue\|_{\frac 53}^{\frac 14}\|\pat \ue\|_2^{\frac {7}{4}} \|\nabla \ue \|_{10}\\ \leq \delta \int_0^\tau\|\pat \nabla \ue \|_{\frac 53}^2  + c_\delta \int_0^\tau \|\pat \ue\|_2^2 \|\nabla \ue\|_{10}^{\frac {8}{7}},
\end{multline*}
for arbitrary $\delta>0$. The first term of the right hand side can be subsumed into the left hand side of \eqref{skorokonec} in virtue of the following estimate

\begin{multline}\label{skorokonec2}
\|\pat\nabla \ue\|_{L^2(L^p)}^2 \le c \|\pat D \ue\|_{L^2(L^p)}^2 = c\int_0^\tau\left(\int_\Omega |\pat D \ue|^p\mu(D\ue)^{\frac p2} \mu(D\ue)^{-\frac p2}\right)^{\frac 2p}\\ \leq c\int_0^\tau\left(\int_\Omega \mu(D\ue) |\pat D \ue|^2\right) \left(\int_\Omega (1+|D\ue|^2)^{\frac p2}\right)^{\frac{2-p}{p}} \\ \leq c\int_{(0,\tau)\times \Omega} \mue(D\ue)|\pat D\ue|^2.
\end{multline}
where Korn's inequality, H\"older's inequality and \eqref{casoveodhady}$_2$ were used. From \eqref{skorokonec} we infer
\begin{multline}\label{casodhad2a}
\int_0^t \pat \|\pat\ue\|_2^2 + c\int_0^t\int_\Omega \mu(|D\ue|) |D\pat \ue|^2\\ \leq c\int_0^t \left(1 + \|\nabla \ue\|_{10}^{\frac{18}{17}}\right) \|\pat \ue \|_2^2 + c \|\pat f\|_{L^2(Q)}^2.
\end{multline}
We know that $\left(1 + \|\nabla \ue\|_{10}^{\frac{8}{7}}\right)\in L^1(0,T)$ as a consequence of \eqref{ieoasgfinal} and embedding $W^{2, \frac 53} \hookrightarrow W^{1,10}$. Thus, we deduce by Gronwall's inequality that
\begin{equation}
\begin{split}\label{casodhad2}
\|\pat \ue\|_{L^\infty(L^2)}&\leq c,\\
\|\pat \nabla \ue \|_{L^2(L^p)}& \leq c.
\end{split}
\end{equation}

\subsection{Consequences}
With $\nabla \ue \in L^{\frac 52}(Q)$ and \eqref{casodhad2}$_2$, we may use the same method as in \eqref{methoda} in order to derive
\begin{equation*}
\|\nabla \ue\|_{L^\infty(L^2)}\leq c.
\end{equation*}

Since we have better regularity of $\nabla \ue$, we may use the same calculations as in \eqref{skorokonec2} in order to derive from \eqref{ieoasg4} and \eqref{casodhad2a} estimates
\begin{equation}
\begin{split}\label{sec.der.u4}
\|\nabla^2 \ue\|_{L^2(L^{\frac4{4-p}})} &\leq c,\\
\|\nabla \pat \ue \|_{L^2(L^{\frac4{4-p}})} &\leq c,
\end{split}
\end{equation}
see also \cite[Section 4, Step 2]{kms3}. We use \eqref{methoda} again and as in \cite[Section 4, Step 3]{kms3} arrive at
\begin{equation}\label{sec.der.u3}
\|\ue\|_{L^\infty(W^{1,p+1})} \leq c.
\end{equation}


\section{Improved time regularity}
\label{imp.tim.reg}

In this section we follow \cite[Section 5, Step 6]{tic}, the same method was used earlier for the problem with periodic boundary conditions in \cite[Section~4, Step~4]{kms3}. To improve information about $\partial_t\ue$ we apply the $L^p$ theory for the classical Stokes system under perfect slip boundary conditions proven in \cite[Section 3]{tic}.

We consider the following system
\begin{equation}\label{gstokes}
\begin{split}
\int_I \langle \partial_t u,\varphi\rangle + \int_Q \mathbb{M}:Du\otimes D\varphi  = \int_Q G: D\varphi\quad \forall \varphi \in L^{q}(I,W^{1,q}_\sigma),
\end{split}
\end{equation}
where the coefficient matrix $\mathbb{M} \in L^{\infty}(Q)$ is symmetric in the sense $M^{kl}_{ij} = M^{ij}_{kl}=M^{ji}_{kl}$ for $i,j,k,l=1,2$ and fulfils for all $B\in\mathbb{R}^{2\times 2}$, $x\in \Omega$ and $t\in I$
$$
\gamma_1|B|^2 \le \mathbb{M}(t,x):B\otimes B \le \gamma_2 |B|^2.
$$
Here the constants $\gamma_1$ and $\gamma_2$ fulfills $0<\gamma_1\le\gamma_2$. In the following lemma that states the $L^q$ theory result, $B^s_{p,q,B,\sigma}$ is the Besov space $B^s_{p,q}$ which is divergence free and reflects the boundary conditions, i.e. $B^s_{q,q,B}=B^s_{q,q}$ provided $s\in(1,1/q)$, $B^s_{q,q,B}=\{u\in B^s_{q,q},\,u\cdot \nu=0 \textrm{ on } \partial\Omega\}$ provided $s\in [1/q,1+1/q)$ and $B^s_{q,q,B}=\{u\in B^s_{q,q},\,u\cdot \nu=0,\, [(Du)\nu]\cdot\tau = 0 \textrm{ on } \partial\Omega\}$ provided $s\in (1+1/q,2]$ and $B^s_{q,q,B,\sigma} = B^s_{q,q,B} \cap L^q_\sigma$. (c.f.  \cite[Corollary 2.6]{steiger}). 

\begin{lemma}{\cite[Lemma 3.11]{tic}}\label{lptheorygs}
Let $\Omega$ be a bounded non-circular $\mathcal{C}^{2,1}$ domain and $q>2$. There exist constants $c_K,c_L>0$ such that if $q\in[2,2+c_L\frac{\gamma_1}{\gamma_2})$, $G\in L^q(L^q)$ and $u_0 \in B^{1-2/q}_{q,q,B,\sigma}$ then the unique weak solution $u\in L^q(I,W^{1,q}_\sigma)$ of \eqref{gstokes} satisfies
\begin{align}
\|Du\|_{L^q(L^q)} + \gamma_2^{-\frac{1}{q}}\|u\|_{BUC(I,B^{1-2/q}_{q,q,B,\sigma})} \le \frac{c_K}{\gamma_1}\Big(\|G\|_{L^q(L^q)} + \gamma_2^{1-\frac{1}{q}}\|u_0\|_{B^{1-2/q}_{q,q,B,\sigma}}\Big).\notag
\end{align}
\end{lemma}

Suppose that $f\in L^{q_1}(W^{-1,q'_1}_\sigma)$ for some $q_1>2$ and $u_0 \in W^{2+\beta, 2}$ for $\beta \in (0,1/4)$. We observe that $(\partial_t \ue,\partial_t \pie)$ solves 

\begin{multline}\label{timederns}
\int_I\langle \partial_t^2\ue,\varphi\rangle +\int_Q \partial_{D\ue}^2\Phi(|D\ue|):D\partial_t \ue\otimes D\varphi \\= \int_I\langle\partial_t (f-(\ue\cdot\nabla)\ue),\varphi\rangle,\quad \forall \varphi\in L^q(I,W^{1,q}_\sigma).
\end{multline}

The right hand side of \eqref{timederns} is bounded uniformly with respect to $\varepsilon \in (0,1)$ in $L^{q_0}(I,W^{-1,{q_0}}_\sigma)$ for some $q_0>2$. Indeed, it follows from the boundedness of $\||\ue| |\partial_t\ue|\|_{L^{q_0}(Q)}$. This holds due to \eqref{sec.der.u3} and interpolation between \eqref{casodhad2}$_1$ and \eqref{sec.der.u4}$_2$.

Set 
$$
V_\varepsilon:=\esssup_Q (1+|D\ue|^2)^{\frac 12}.
$$ 
The growth conditions \eqref{trunc.growth.con} imply
$$
cV_\varepsilon^{p-2}(|A|)|B|^2 \le \partial^2_A\Phie(|A|):B\otimes B \le \gamma_2 |B|^2.
$$

Lemma \ref{lptheorygs} with $\gamma_1= cV_\varepsilon^{p-2}(|A|)$ provides the existence of positive constants $c_K$ and $c_L$  such that for all $q\in (2,q_2]$, where $q_2:= 2+c_LV_\varepsilon^{p-2}$ holds
\begin{multline}\label{est.v1}
\|\nabla \partial_t \ue\|_{L^{q}(Q)} + \gamma_2^{-\frac{1}{q}}\|\partial_t \ue\|_{BUC(I,B^{1-2/q}_{q,q,B,\sigma})} \\ \le \frac{c_K}{c_5}V_\varepsilon^{2-p}\Big(\|f\|_{L^q(I,W^{-1,q'})} + \gamma_2^{1-\frac 1q}\|\partial_tu_0\|_{B^{1-2/q}_{q,q,B,\sigma}}\Big).
\end{multline}
Without loss of generality we may assume that $q_2<q_0$. It remains to show that last norm on the right hand side of \eqref{est.v1} is bounded. This can by done in the same way as in \cite[Section 4, Step 3]{tic}, therefore we omit it.

As $V_\varepsilon^{p-2}\leq 1$, we may assume without loss of generality that $q_2 = 2+ c_LV_\varepsilon^{p-2}\leq 3$ for some suitable $c_L$. The estimate \eqref{est.v1} therefore reduces to

\begin{equation}\label{cas.derivace.odhad}
\|\partial_t\ue\|_{BUC(I,B^{1-2/q}_{q,q,B,\sigma})} \le cV_\varepsilon^{2-p}.
\end{equation}
By interpolation we get from \eqref{sec.der.u3} and \eqref{cas.derivace.odhad}
\begin{equation}\label{interpolace.derivace}
\|\partial_t u^\varepsilon\|_{L^\infty(L^{q_1})} \leq c\left(1+V_\varepsilon^{b(2-p)}\right), \quad \mbox{where }b = \frac {q_2}{q_1}\frac{q_1 - 2}{q_2-2}
\end{equation}
for every $q_1 \in [2,q_2]$.

\section{Improved space regularity}
\label{imp.spa.reg}

In order to provide an estimate of the second gradient of $u^\varepsilon$, keeping at our disposal \eqref{interpolace.derivace} we move $\partial_t \ue$ to the right hand side and use a stationary theory for a quadratic growth with lower bound expressed in terms of $V_\varepsilon$. 

In general, we investigate the problem
\begin{equation}\label{stacgstokes}
\int_\Omega\mathbb{M}:Du\otimes D\varphi = \int_\Omega G:D\varphi \quad \forall \varphi \in W_\sigma^{1,q},
\end{equation}
where $q\ge 2$, coefficient matrix $\mathbb{M} \in L^\infty$ is symmetric and fulfills for all $B\in\mathbb{R}^{2\times 2}$, $x\in \Omega$ and $0<\gamma_1\le \gamma_2$ the growth
$$
\gamma_1|B|^2 \le \mathbb{M}(x):B\otimes B \le \gamma_2 |B|^2.
$$

We will need
\begin{lemma}{\cite[Lemma 3.12]{tic}}\label{staclptheory}
Let $\Omega$ be a bounded non-axisymmetric $\mathcal{C}^{2,1}$ domain. Then there are constants $c_K,c_L>0$ such that if $q\in[2,2+c_L\frac{\gamma_1}{\gamma_2})$ and $G \in L^{q}$, then the unique weak solution of \eqref{stacgstokes} satisfies
$$
\|D u\|_{L^q} \le \frac{c_K}{\gamma_1}\|G\|_{L^q}.
$$
\end{lemma}

As we would follow line by line the process in \cite[Step~7 in Section~5]{tic}, we restrict ourselves only on commentary of the important steps. In spite of \cite{tic}, where the author is involved in the super-quadratic case and therefore $\gamma_2$ is dependent on the truncation parameter $\varepsilon$, we are interested in the precise dependence of $\gamma_1$. Thus, we focus on the parts of the proof where the constants $\gamma_1$ and $\gamma_2$ plays the key role.

To obtain boundary regularity result we localize the problem and work in $\Omega_{3r}$. At first we gain information in a tangent direction. After we subsume $\partial_t\ue$ into the right hand side, we test the equation $\eqref{ap.sys}_1$ by the tangent derivative of the test function, precisely by $-\pt \varphi \xi$ with $\varphi \in W^{1,q'}_\sigma$, $\supp \varphi \subset \overline{\Omega_{3r}}$. After some manipulation involving integration by parts and Bogovski\u {\i} type correction we arrive at

$$
\int_{\Omega_{3r}} \partial_{Du^\varepsilon}\ese(Du^\varepsilon):D\pt\ue\xi\otimes D\varphi  = \langle g,\varphi \rangle, \qquad \forall \varphi \in W^{1,q'}_{\sigma}
$$
with
$$
\|g\|_{W^{-1, q_3}_\sigma}\leq c\left(1+ V_\varepsilon^{(2-p)b}\right)
$$ 
for some $q_3\in (2,q_2)$. Application of Lemma \ref{staclptheory} yields
\begin{equation}\label{tec.druhy.grad}
\|\nabla \pt u^\varepsilon \xi\|_{L^{q_3}} \leq c V_\varepsilon^{2-p}(1+V^{(2-p)b}_\varepsilon) + c.
\end{equation}
Moreover, $q_3$ can be expressed as a convex combination of $q_2$ and $2$, i.e. there exists $\lambda_1\in (0,1)$ such that for every $\lambda \in (0,\lambda_1)$ \eqref{tec.druhy.grad} holds provided
$$
q_3 = \lambda q_2 + (1-\lambda) 2.
$$
In this case 
\begin{equation}\label{rovnice.na.a}
b= 1- \frac{(1-\lambda)2}{\lambda q_2 + (1-\lambda)2}\leq 1-\frac{(1-\lambda)2}{\lambda 3 + (1-\lambda)2}
\end{equation}
and, due to expression of $q_2$, 
\begin{equation}\label{gamma.1}
1 - \frac 2{q_3} \geq c\lambda V_\varepsilon^{p-2}
\end{equation}
with some $c$ independent of $\varepsilon$.

To control whole second gradient of $\ue$ it remains to obtain estimate of the type \eqref{tec.druhy.grad} in the normal direction. As it is possible to extract $\partial_2^2 \ue_2$ from $\diver \ue=0$, it suffices to control $\partial_2^2 \ue_1$. In the same way as in \cite[Theorem~3.19]{kms2} or \cite[Step~7 of Section~5]{tic} we would apply the operator $\curl$ to $\eqref{ap.sys}_1$ and use growth conditions and Ne\v cas' theorem on negative norms and $\partial_{12}\ese_{12}\geq V_\varepsilon^{p-2}$.

Finally, after passing from the local estimates to the whole domain $\Omega$ we arrive at

\begin{equation}\label{odhad.na.druhou.der}
\|\nabla^2 u^\varepsilon\|_{L^{q_3}}\leq c(1+V_\varepsilon^{2-p})\left(1 + V_\varepsilon^{(2-p)(1+b)}\right)\leq c \left(1+ V_\varepsilon^{(2-p)(2+b)}\right)
\end{equation}
for almost all $t\in I$.

\section{Proof of the main theorem}

We use lemma, whose proof follows from \cite[proof of Theorem 2.4.1]{ziemer}.
\begin{lemma}
Let $\Omega\subset \mathbb{R}^2$ be a bounded $\mathcal{C}^2$ domain and $v\in W^{1,q}$ for some $q>2$. Then $v\in \mathcal{C}(\overline{\Omega})$ and there is $c>0$ independent of $q$ such that
\begin{equation*}
\sup_{\Omega}|v| \le c\left(\frac{1}{1-\frac 2q}\right)^{1-1/q}\|v\|_{1,q}.
\end{equation*}
\end{lemma}

With the help of above mentioned lemma, \eqref{gamma.1} and \eqref{odhad.na.druhou.der} we derive
\begin{multline}\label{lnekonecno}
V_\varepsilon \leq c\left(\frac 1{1-\frac 2{q_3}}\right)^{1-\frac 1{q_3}} \esssup_t \|\nabla^2 u^\varepsilon(t)\|_{q_3} + c\\ \leq c V_\varepsilon^{(2-p)(1-\frac 1{q_3})} \left(1 + V_\varepsilon^{(2-p)(2 + b)}\right) + c.
\end{multline}
As $p\geq \frac 53$ and $q_3<3$, $\lambda\in (0,1)$ can be chosen such that $b<\frac 13$ (see \eqref{rovnice.na.a}) and, consequently, 
$$
(2-p)\left(3-\frac 1{q_2} + b\right)<1.
$$
Further, \eqref{lnekonecno} yields $V_\varepsilon < c$ with $c$ independent of $\varepsilon$. Thus, for $\varepsilon$ sufficiently small, a solution $\ue$ to \eqref{ap.sys} coincides with a solution $u$ to \eqref{ns} and satisfies
$$
u\in L^\infty (I,W^{2,q})\quad \mbox{and}\quad \partial_t u\in L^\infty(I, W^{1,q})
$$
for some $q>2$ and we obtain $u\in \mathcal C^{1,\alpha}(Q)$ for some $\alpha>0$ in virtue of  

\begin{lemma}{\cite[Lemma 2.2]{JoSta}}
Assume that for any $\beta>0$ and $r>1$ $f\in L^{\infty}(\mathcal C^{1,\alpha})$ and  $\partial_t f \in L^r(W^{1,r})$. Then for $\alpha = \min\{\beta, \beta(r-1)/(\beta r + 2)\}$
$$
f \in \mathcal C^{1,\alpha}(Q).
$$
\end{lemma}

Since the bound on $V_\varepsilon$ allows to pass from sub-quadratic to quadratic case, the regularity of pressure can be deduced similarly as in \cite[Section~4]{tic}.


\end{document}